\documentclass[12pt]{amsart}
\usepackage{}
\usepackage{amsmath}
\usepackage{amsthm}
\usepackage{amsfonts,amssymb}
 \usepackage{color}
\usepackage{enumerate}
\RequirePackage{srcltx}

\newtheorem{thm}{Theorem}[section]
\newtheorem{cor}[thm]{Corollary}
\newtheorem{lem}[thm]{Lemma}

\newtheorem{defn}[thm]{Definition}
\theoremstyle{remark}
\newtheorem{rem}{Remark}[section]

\newcommand{\thmref}[1]{Theorem~\ref{#1}}

\numberwithin{equation}{section}

 \def\tr{{\triangle}}
 
\def\sub{\substack}
\def\sph{\mathbb{S}^{d-1}}
\def\f{\frac}

\def\Bl{\Bigl} \def\Br{\Bigr}

 \def\bl{\bigl}
\def\br{\bigr}
\def\({\left(}
\def \){ \right)}
\def\sa{\sigma}
 \def\a{{\alpha}}
 \def\b{{\beta}}

 \def\t{{\theta}}
 \def\l{{\lambda}}
 \def\ld{\lambda}

 \def\s{{\sigma}}

 \def\CH{{\mathcal H}}

 \def\NN{{\mathbb N}}

 \def\RR{{\mathbb R}}
  \def\R{{\mathbb R}}
  \def\SS{{\mathbb S}}
 
 \def\al{\a}
 
 \def\proj{\operatorname{proj}}

\def\hb{\hfill$\Box$}

\def\leqs{\leqslant}

\newcommand{\wh}{\widehat}

\def\sph{\mathbb{S}^{d-1}}

\def\s{\sa}

\def\Ga{\Gamma}

\def\be{\b}
\def\HH{\mathcal{H}}
\begin{document}

\title[Reverse H\"{o}lder's inequality]
{
Reverse H\"{o}lder's  inequality for spherical harmonics
}
\author{Feng Dai}
\address{Department of Mathematical and Statistical Sciences\\
University of Alberta\\ Edmonton, Alberta T6G 2G1, Canada.}
\email{fdai@ualberta.ca}
\author{Han Feng}
\address{Department of Mathematical and Statistical Sciences\\
University of Alberta\\ Edmonton, Alberta T6G 2G1, Canada.}
\email{hfeng3@ualberta.ca}
\author{Sergey Tikhonov} \address{ICREA, Centre de Recerca Matem\`{a}tica\\
Campus de Bellaterra, Edifici C
08193 Bellaterra (Barcelona), Spain, and Universitat Aut\`{o}noma de Barcelona.}
\email{ stikhonov@crm.cat}

\date{\today}
\keywords{spherical harmonics, polynomial inequalities, restriction theorems} \subjclass{33C50,
33C52, 42B15, 42C10}

\thanks{The first and the second   authors  were partially supported  by the NSERC Canada
under grant RGPIN 311678-2010. The third author was partially supported  by
 MTM 2011-27637, 2014 SGR 289,  RFFI 13-01-00043 and
the Alexander von Humboldt Foundation.
}
\begin{abstract}
This paper determines the sharp asymptotic order   of
 the following  reverse H\"{o}lder  inequality for  spherical harmonics $Y_n$  of degree $n$  on the unit sphere $\sph$ of $\RR^d$ as $n\to \infty$:
 $$\|Y_n\|_{L^q(\sph)}\leq C n^{\a(p,q)}\|Y_n\|_{L^p(\sph)},\   \
   0<p<q\leq \infty. $$   In many cases, these  sharp estimates turn out to be  significantly  better   than  the corresponding estimates in the  Nilkolskii inequality  for spherical polynomials.  Furthermore, they    allow us to improve two recent results on the restriction conjecture and the
    sharp Pitt  inequalities  for the Fourier transform  on $\RR^d$.
\end{abstract}

\maketitle
\smallskip

\section{Introduction}

Let  $\sph=\{ x\in \mathbb{R}^{d}:   \|x\|=1\}$ denote the
unit sphere of $\RR^d$ endowed  with the usual Haar measure
  $d\sa(x)$, where $\|\cdot\|$ denotes the Euclidean norm of $\RR^d$.     Given $0<p\leq \infty$, we denote by $L^p(\sph)$ the usual Lebesgue $L^p$-space defined with respect to the measure $d\s(x)$ on $\sph$, and by $\|\cdot\|_p$ the norm of $L^p(\sph)$.
Throughout the paper, unless otherwise stated, all functions on $\sph$ will be assumed to be real-valued and measurable, and   the
notation $ A\sim B$  means  that there exists an inessential
constant $c>0$, called the constant of equivalence, such that
$ c^{-1} A \leq B \leq c A.$

Let $ \Pi_n^d$  denote the space of all
spherical polynomials of degree at most $n$ on $\sph$  (i.e., restrictions  on  $\sph$ of  polynomials  in $d$ variables of total
degree at most $n$), and  $\HH_n^d$  the space of all spherical harmonics of degree $n$ on $\sph$.  As is well known (see, for instance, \cite[chapter 1]{BOOK}),   $\HH_n^d$ and $\Pi_n^d$ are all  finite dimensional spaces with $\text{dim} \HH_n^d\sim n^{d-2}$ and $\text{dim} \Pi_n^d \sim n^{d-1}$ as $n\to \infty$. Furthermore, the spaces $\HH_k^d$, $k=0,1,\cdots$ are mutually orthogonal with respect to the inner product of $L^2(\sph)$, and
each space $\Pi_n^d$ can be written as a direct sum $\Pi_n^d =\sum_{j=0}^n \HH_j^d$.
Since the space of spherical polynomials is dense in $L^2(\sph)$,    each $f\in L^2(\sph)$ has a spherical harmonic  expansion,
$
f=\sum_{k=0}^\infty \proj_k f,$
   where $\proj_k$ is the orthogonal projection of $L^2(\sph)$ onto the space $\HH_k^d$ of spherical harmonics.  The orthogonal projection $\proj_k$   has an integral representation:
   \begin{equation}\label{1-1-0}
    \proj_kf(x) =C_{k,d}\int_{\sph} f(y) P_k^{(\f {d-3}2, \f
{d-3}2)}(x\cdot y)\, d\sa(y),\ \  \ x\in\sph,
\end{equation}
where
$$C_{k,d}:= \f {\Ga(\f d2)\Ga(\f {d-1}2) }{2\pi^{d/2}\Ga(d-1) } \f { (2k+d-2)\Ga(k+d-2)}{\Ga(k+\f {d-1}2)},$$
and    $P_k^{(\al,\be)}$ denotes   the usual Jacobi polynomial of degree $k$ and   indices
$\al, \be$, as defined in \cite[Chapter IV]{Sz}.

Our goal  in this paper is to find  a sharp asymptotic order  of the quantity  $\sup_{Y_n\in\HH_n^d}\f{ \|Y_n\|_q}{\|Y_n\|_p}$  for  $0<p<q\leq \infty$ as $n\to \infty$. The background of this problem is as follows.  In 1986,  Sogge \cite{So1} proved that for $d\ge 3$, and $\l:=\f {d-2}2$,
\begin{equation}\label{1-1}
\sup_{Y_n\in\mathcal{H}_n^d}\f{\|Y_n\|_{L^q(\sph)}}{\|Y_n\|_{L^2(\sph)}} \sim   \begin{cases}
  n^{\l(\f 12-\f 1q)}, &  \qquad \quad 2\leq q\le  2(1+\f 1\l),\\
  n ^{2\lambda (\f 12-\f 1q)-\f1q}, & \qquad\quad  2(1+\f 1\l)\le q\le \infty,\end{cases}
\end{equation}
which  confirms  a conjecture of Stanton--Weinstein \cite{Stan} in the  case of  $d=3$ and $q=4$.  Here and throughout the paper, it is agreed that $0/0=0$.
 Recently, De Carli and Grafakos  \cite{grafakos}  proved that if $1\leq p\leq q\leq 2$ and  $Y_n\in \HH_n^d$ can be written in the form
\begin{equation}\label{1-3}Y_n(x) = e^{ i m_{d-2} x_{d-1}} \prod_{k=0}^{d-2}
(\sin x_{k+1})^{m_{k+1}} P_{m_k-m_{k+1}}^{(m_{k+1}+\f {d-2-k}2, m_{k+1}+\f {d-2-k}2)} (\cos x_{k+1}),\end{equation}
with  $n=m_0 \ge m_1\ge \cdots m_{d-2}\ge 0$  being integers, then
 \begin{equation}\label{1-4}
 \frac{ \|Y_n\|_{L^q(\mathbb{S}^{d-1})}  }{ \|Y_n\|_{L^p(\mathbb{S}^{d-1})}  }\leqs C n^{\frac{d-2}{2} (\f {1}{p}-\f 1q)},\    \    \  \   1\leq p<q\leq 2,  \end{equation}
 which  was further  applied  in  \cite{grafakos}  to prove  the restriction conjecture for the class of functions consisting of
products of radial functions and spherical harmonics that are  in the form \eqref{1-3}.
Note  that  the set of functions  $Y_n$ in \eqref{1-3} with $n=m_0 \ge m_1\ge \cdots m_{d-2}\ge 0$  forms a linear basis of the space $\HH_n^d$.
It is therefore natural to ask whether or not  \eqref{1-4} holds for all spherical harmonics $Y_n$ of degree $n$.    A related work in this direction was  done recently by De Carli, Gorbachev  and  Tikhonov  in  \cite{DeC},
 where  the following weaker estimate  was obtained  for all spherical harmonics and   applied   to  study  a sharp Pitt inequality for the Fourier transform on $\RR^d$:
\begin{equation}\label{1-2}
\sup_{Y_n\in\HH_n^d} \f { \|Y_n\|_{p'}} { \|Y_n\|_p} \leq C n^{ (d-1) (\f 1p-\f 12)},\    \   \    \   \f1p +\f 1{p'}=1,\    \    \  1\leq p\leq 2,
\end{equation}
 Finally,   let us recall the following  well-known result of  Kamzolov  \cite{Kam}  on the  Nikolskii inequality  for  spherical polynomials:
  \begin{equation}\label{Nikol}
 \|P_n\|_q \leqs C n^{(d-1) (\f 1p-\f1q)}\|P_n\|_p,\qquad \forall  P_n\in\Pi_n^d,\   \     0<p<q\leq \infty.
 \end{equation}
 Since $\HH_n^d\subset \Pi_n^d$, the Nikolskii inequality \eqref{Nikol} is applicable to every spherical harmonics $Y_n\in \HH_n^d$. It  turns out, however, that the resulting estimates are not sharp for spherical harmonics in many cases (see, for instance, \eqref{1-1}, \eqref{1-2} and \eqref{1-4}).

 In this paper, we will prove  the following result, which, in particular, shows that \eqref{1-4} holds for all spherical harmonics $Y_n\in\HH_n^d$, and the upper bound on the right hand side of   \eqref{1-2} can be improved to be
$ C n^{ (d-2) (\f 1p-\f 12)}$.

\newpage\begin{thm}\label{cor2} Assume that
  $d\ge 3$ and  $\f 1p+\f 1{p'}=1$  if $p\ge 1$.  Set  $\l:=\f {d-2}2$.
\begin{enumerate}[\rm (i)]
\item
If  either  $0<p\leq 1$ and $p<q\leq \infty$, or $1\leq  p\leq 2$ and $ p<q \leq
 (1+\f 1\l)p'$,
then
\begin{equation}\label{1-6}
\sup_{\sub{Y_n\in\mathcal{H}_n^d}}\f{\|Y_n\|_q}{\|Y_n\|_p} \sim   n^{\l(\f 1p-\f 1q)}.\end{equation}

\item If either $1\leq p\leq 2$ and $q\ge (1+\f 1\l)p'$, or
$2\leq p< 2+\f 1\l$ and  $q> 2+\f 2\l$, then
$$\sup_{\sub{Y_n\in\mathcal{H}_n^d}}\f{\|Y_n\|_q}{\|Y_n\|_p}\sim n ^{2\lambda (\f 12-\f 1q)-\f1q}.$$

\item If $ 2+\f 1\l<p <q \leq \infty$, then
$$\sup_{\sub{Y_n\in\mathcal{H}_n^d}}\f{\|Y_n\|_q}{\|Y_n\|_p}\sim  n^{(2\l+1) (\f 1p-\f1q)}.$$

\item  If  $d=3$  and  $2\leq p< 4=2+\f1\l$, then for  $q\ge 3 p'=(1+\f 1\l) p'$,
$$\sup_{\sub{Y_n\in\mathcal{H}_n^d}}\f{\|Y_n\|_q}{\|Y_n\|_p}\sim
n^{\f 12 -\f 2q},
$$
whereas for $p<q\leq 3p'$,
$$\sup_{\sub{Y_n\in\mathcal{H}_n^d}}\f{\|Y_n\|_q}{\|Y_n\|_p}
 \sim  n^{\f 12 (\f 1p-\f1q)}.$$

\end{enumerate}
\end{thm}

Of particular interest is the case when $1\leq p\leq 2$ and $q=p'$, where our result can be stated as follows:

\begin{cor}\label{cor-1-2}If $Y_n\in \mathcal{H}_n^d$ and $1\leq p\leq 2$, then
\begin{equation}\label{1-7} \|Y_n\|_{p'}\leq C n^{(d-2)(\f 1p-\f12)} \|Y_n\|_p,\   \   \  1\leq p\leq 2.\end{equation}
Furthermore, this estimate is sharp.
\end{cor}

Several remarks are in order.

 \begin{rem} Estimate  \eqref{1-7}   for $p=p_\l:=1+\f{\l}{\l+2}$ follows directly from the
  well-known result of Sogge \cite{So1} on  the orthogonal projection
   $\proj_n : L^2(\sph)\to \HH_n^d$.  However, for $1\leq p<2$ and $p\neq p_\l$,
   the sharp estimate   \eqref{1-7} in Corollary \ref{cor-1-2}  is nontrivial and cannot be deduced from the result of Sogge \cite{So1}.
 Indeed, it was shown in \cite{So1}
 that for $1\leq p\leq p_\l:=1+\f{\l}{\l+2}$,
\begin{equation}\label{lemma-sogge}
\|\proj_n f\|_{2}\leq C n^{\l(\f{1}{p}-\f{1}{2} )+ \f{1}{2p(\l+2)}(p_\l-p)}
 \|f\|_{p},\    \   \   \    \forall  f\in L^{p}(\sph),
 \end{equation}
and this estimate is sharp.  Since $\proj_n f=f$ for $f\in \HH_n^d$,  this leads to the inequality
 $$\|Y_n\|_{2} \leq C n^{\l(\f{1}{p}-\f{1}{2} )+ \f{1}{2p(\l+2)}(p_\l-p)}\|Y_n\|_{p},\   \   \    \    \forall\, Y_n\in\HH_n^d,\    \   1\leq p\leq p_\l,$$
 which, according to Corollary \ref{cor-1-2},  is not sharp unless $p=p_\l$.

 \end{rem}
 \begin{rem}
  Interesting reverse H\"{o}lder inequalities for spherical harmonics,
$$
\sup_{Y_n\in\HH_n^d} \f{\|Y_n\|_q}{\|Y_n\|_p} \leq C(n,q)
$$
with the constant  $C(n,q)$ being independent of the dimension $d$ but dependent on the degree $n$ of spherical harmonics,  were obtained in \cite{Duo} for some pairs of $(p,q)$, $0<p<q<\infty$.  The general constants $C$ in our paper are dependent on the dimension $d$, but independent of the degree $n$.

\end{rem}

 \begin{rem}    For $d\ge 4$,  it remains open  to find the
 asymptotic estimate of the supremum on the left hand side of \eqref{1-6} for  $2<p<1+\f 1\l$ and $p<q< 2+\f 2\l$.
\end{rem}

This  paper is organized as follows.
In Section 2, we construct a sequence of convolution  operators $\{T_n\}_{n=0}^\infty$ on $L^1(\sph)$ with the properties that
    $T_n f=f$ for $f\in\HH_n^d$,    $|T_n f|\leq C \sup_{0\leq j \leq d} |\proj_{n+2j} f|$
    and
    $\|T_n f\|_\infty \leq C n^\l \|f\|_1$  for all    $f\in L^1(\sph)$.
These operators  play an indispensable role in the proof of  Theorem \ref{cor2},
which is given in the third section.
Finally,  in Section 4,   we give two applications  of  our main result,   improving  a recent result of \cite{grafakos} on restriction conjecture and a result of \cite{DeC} on sharp Pitt's inequality.

\section{ A sequence of convolution operators  }

 We start with the following well-known result  of Sogge \cite{So1}   on  the operator norms  of the orthogonal projections
  $\proj_n:  L^2(\sph)\to \HH_n^d$.

  \begin{lem} \cite{So1}\label{lem-4-3}
Let $n\in\NN$, $d\ge 3$ and $\l=\f{d-2}2$. Then the following statements hold: \begin{enumerate}[\rm (i)]

\item   If  $1\leq p\leq p_\l:=1+\f \l {\l+2}$, then
$$\|\proj_n f\|_2\leq C n^{(2\l+1)(\f 1p-\f12)-\f12}\|f\|_p.$$
\item   If  $p_\l\leq p\leq 2$, then
$$\|\proj_n f\|_2\leq C n^{\l(\f 1p-\f12)}\|f\|_p.$$

\item  If  $2+\f 2\l\leq q \leq \infty$, then
$$\|\proj_n f\|_q\leq C n^{(2\l+1)(\f12-\f1q)-\f12}\|f\|_2.$$

\item If $2\leq q \leq 2+\f 2\l$, then
$$\|\proj_n f\|_q \leq C n^{\l(\f 12-\f 1q)} \|f\|_2.$$
Here, the letter $C$ denotes a general positive constant  independent of  $n$ and $f$.
\end{enumerate}
\end{lem}

  As was pointed out  in the introduction,  Lemma \ref{lem-4-3}   will  not  be enough for the proof of our main result.
   The crucial step in the proof of Theorem \ref{cor2} is  to  construct a sequence of linear operators $\{T_n\}_{n=0}^\infty$ with the properties that
    $T_n f=f$ for $f\in\HH_n^d$,    $|T_n f|\leq C \sup_{0\leq j \leq d} |\proj_{n+2j} f|$  and
    $\|T_n f\|_\infty \leq C n^\l \|f\|_1$  for all    $f\in L^1(\sph)$ .

    To define  the operators   $T_n$,  we need  to recall  several   notations.
First, given  $h\in\NN$,
and
a sequence $\{a_n\}_{n=0}^\infty$ of real numbers,  define
$$\tr_h a_n =a_n-a_{n+h},\    \
\tr_h^{\ell+1}=\tr_h \tr_h^\ell,\   \  \ell=1,2,\ldots.$$
Next, let $$R_n^\l (\cos\t):=\f {P_n^{(\l-\f12, \l-\f12)}(\cos\t) }{P_n^{(\l-\f12, \l-\f12)}(1)},\    \   \   \t\in [0,\pi]$$
denote the normalized Jacobi polynomial, and  for   a step $h\in \NN$, define
$$\tr_h^\ell R_n^\l(\cos\t)=\tr_h^\ell a_n,\    \   \ell=1,2,\ldots,\    \      n=0,1,\cdots, $$
with  $a_n:= R_n^{\l}(\cos\t)$.  Here and throughout, the difference operator in $\tr_h^\ell R_n^\l(\cos\t)$ is always acting on the integer $n$.
In the case when the step $h=1$, we have
the following estimate  (\cite[Lemma B.5.1]{BOOK}, \cite{DD}):
 \begin{equation}\label{0-1}
 \Bl|\tr_1^{\ell}  R_n^\l (\cos\t)
 \Br|\leq C \t^\ell (1+n\t)^{-\ld},\   \    \  \t\in [0, \pi/2],\   \ \ell\in \NN.
 \end{equation}
On the other hand, however, the $\ell$-th order
difference $\tr_1^{\ell}  R_n^\l (\cos\t)$  with step $h=1$  does not provide a desirable  upper  estimate when
    $\t$ is  close to $\pi$, and as will be seen in our later proof,
 estimate \eqref{0-1} itself will not be enough for our purpose.

To overcome this difficulty, instead of the difference with step $1$,
 we consider the $\ell$-th order  difference $\tr_2^{\ell}  R_n^\l (\cos\t)$ with step $h=2$.
     Since
 $\tr_2^\ell a_n =\sum_{j=0}^\ell\binom{\ell} j \tr_1^\ell a_{n+j},$
 on one hand, \eqref{0-1}   implies that
 $$ \Bl|\tr_2^{\ell}  R_n^\l (\cos\t)
 \Br|\leq C \t^\ell (1+n\t)^{-\ld},\    \  \t\in [0,\pi/2].$$
On the other hand,
 however, since
 $$\tr_2^{\ell}  R_n^\l (\cos\t) =\sum_{j=0}^{\ell} (-1)^j \binom{\ell} j R_{n+2j}^\l (\cos \t),$$
  and since $R_{n+2j}^\l (-z) = (-1)^n R_{n+2j}^\l (z)$,
  we have
  $\tr_2^{\ell}  R_n^\l (\cos(\pi -\t))=(-1)^n \tr_2^{\ell}  R_n^\l (\cos\t).$
 It follows that
 \begin{equation}\label{0-2}
\Bl|\tr_2^{\ell}  R_n^\l (\cos\t)
 \Br|\leq C \begin{cases}
 \t^\ell (1+n\t)^{-\ld},\    \  \t\in [0, \pi/2],\\
 (\pi-\t)^{\ell} (1+n (\pi-\t))^{-\ld},\    \  \t\in [\pi/2, \pi].\end{cases}
 \end{equation}

By \eqref{1-1-0}, we obtain   that  for every $P\in \HH_n^d$,
$$ P(x)=c_n \int_{\sph} P(y) R_n^\l (x\cdot y) \, d\s(y),\   \  x\in\sph,
$$
where
$$c_n:=\f{\Ga(\f d2)}{2\pi^{d/2}}\f{ d+2n-2}{d+n-2}\f{\Ga( d+n-1)}{\Ga(n+1)\Ga(d-1)}\sim n^{d-2},$$
and $x\cdot y$ denotes the dot product of $x, y\in\RR^d$.
Since
$R_j^\l (x\cdot)\in\HH_j^d$ for any fixed $x\in\sph$,
it  follows by the orthogonality of spherical harmonics
that for any
 $P\in \mathcal{H}_n^d$, and   any $\ell\in\NN$,
 \begin{align}
   P(x)&
   = c_n \sum_{j=0}^\ell (-1)^j \binom {\ell}j \int_{\sph} P(y) R_{n+2j}^\l (x\cdot y)\, d\s(y)\notag \\
   &=c_n \int_{\sph} P(y) \tr_2^\ell
   R_n^\l (x\cdot y) \, d\s(y).\label{2-1}
 \end{align}
 For the rest of the paper, we will choose $\ell$ to be an integer bigger than $\l$ (for instance, we may  set $\ell=d-2$),
  so that by \eqref{0-2}, we have
  \begin{equation}\label{2-4-1}
    \Bl|\tr_2^{\ell}  R_n^\l (\cos\t)
 \Br|\leq C n^{-\l}.
  \end{equation}

Now we are in a position to define the operators $T_n$.
 \begin{defn}For $f\in L(\sph)$, we define
 \begin{equation}\label{2-2}
 T_n f(x):=\int_{\sph} f(y) \Phi_n(x\cdot y)\, d\s(y),\
 \ x\in\sph,\end{equation}
 where
 $$\Phi_n(\cos \t):=
 c_n \sum_{j=0}^{d-2} (-1)^j \binom {d-2}jR_{n+2j}^\l(\cos\t).$$
 \end{defn}

By \eqref{2-4-1},  we have
 \begin{equation}\label{2-6-1}
 |\Phi_n(\cos\t)|\leq C n^\l,\   \    \
 \t\in [0,\pi],\end{equation}
 whereas by \eqref{2-1}
 \begin{equation}\label{2-7-1}
    T_n P(x) =P(x),\   \   \  \forall P\in \HH_n^d, \   \   \forall x\in \sph.
 \end{equation}

The main result of  this section can now be stated as follows.

\begin{thm}\label{thm-2-3}
\begin{enumerate}[\rm (i)]
\item If  $1\leq p\leq 2$ and $p'\leq q\leq (1+\f 1\l) p'$, then
 \begin{equation}\label{2-5}
    \|T_nf\|_{q}\leq C n^{\l(\f 1p-\f1q)} \|f\|_p,\   \  \forall f\in L^p(\sph).
 \end{equation}
 \item If $1\leq p\leq 2$ and $q\ge (1+\f 1\l) p'$, then
 \begin{align*}
    \|T_n f\|_q \leq C n^{\l-\f{2\l+1}q}\|f\|_p,\   \   \forall f\in L^p(\sph).
 \end{align*}
  \end{enumerate}
\end{thm}
\begin{proof} First, we prove the assertion (i).
Note   that  by definition,
 for each $f\in L^2(\sph)$,
\begin{equation}\label{3-4-0}T_n f =
 \sum_{j=0}^{d-2} (-1)^j \binom {d-2}j \f{c_n}{c_{n+2j}}
 \proj_{n+2j} f,\end{equation}
 which  implies that
 \begin{equation}\label{2-4}
    \|T_n f\|_2\leq C \|f\|_2,\   \  \forall f\in L^2(\sph).
 \end{equation}
 On the other hand, however, using \eqref{2-6-1}, we have
\begin{equation}\label{2-3}
    \|T_n f\|_\infty\leq C n^\l\|f\|_1,\
     \  \forall f\in L^1(\sph).
\end{equation}
 Thus, applying the Riesz-Thorin interpolation theorem,
 and using \eqref{2-4} and \eqref{2-3},
     we deduce that for $1\leq p\leq 2$,
 \begin{equation}\label{2-5-0}
    \|T_nf\|_{p'}\leq C n^{(d-2)(\f 1p-\f12)} \|f\|_p,\   \  \forall f\in L^p(\sph).
 \end{equation}

 Next, by (iv) of Lemma \ref{lem-4-3}, and
 using \eqref{3-4-0}, we obtain  that for $2\leq r\leq 2(1+\f 1\l)$,
 \begin{equation}\label{2-13}
    \|T_n f\|_{r} \leq C n^{\l( \f 12-\f 1{r})}\|f\|_2,\   \   \forall f\in L^2(\sph).
 \end{equation}
 Assume that $1\leq p\leq 2$ and $p'\leq q\leq (1+\f 1\l) p'$.
 Let $\t=\f 2{p'}\in [0,1]$, and let $r=\t q=\f 2{p'} q$.
 Then $2\leq r\leq 2(1+\f 1\l)$, and
 $$\f 1p = 1-\t +\f \t 2,\     \    \f 1q =\f{1-\t}{\infty}+\f \t r.$$
 Thus, by \eqref{2-5-0}, \eqref{2-13}
     and applying the Riesz-Thorin interpolation theorem, we obtain that
 $$\|T_n f\|_q \leq
  C n^{\l (1-\t)} n^{\l (\f 12-\f 1r)\t}\|f\|_p= C n^{\l (\f 1p-\f 1q)}\|f\|_p.$$
This completes the proof of the assertion (i).

Assertion (ii) can be proved similarly.
Indeed, using \eqref{3-4-0} and  (iii)
of Lemma \ref{lem-4-3},
 we have that for
$r\ge  2(1+\f 1\l)$,
 \begin{equation}\label{2-14}
    \|T_n f\|_{r} \leq C n^{2\l (\f 12 -\f 1r)-\f 1r}\|f\|_2,\   \   \forall f\in L^2(\sph).
 \end{equation}
 Assume that $1\leq p\leq 2$ and $ q\ge (1+\f 1\l) p'$.
 Let $\t=\f 2{p'}$ and $r=\t q =\f {2}{p'} q$. Then
 $r\ge 2(1+\f 1\l)$. Using \eqref{2-14}, \eqref{2-5-0}
  and
  applying  the Riesz-Thorin interpolation theorem,
  we deduce that
      \begin{align*}
    \|T_n f\|_q &\leq C n^{\l(1-\t)} n ^{ (d-2)\t (\f 12 -\f 1r)-\f \t r} \|f\|_p=C n^{\l-\f{2\l+1}q}\|f\|_p\\
    &=C n^{(d-2) (\f 12-\f1q)-\f1q}\|f\|_p.
 \end{align*}
This completes the proof of (ii).

\end{proof}

\section{Proof of Theorem  \ref{cor2}}

 The stated lower estimates of Theorem \ref{cor2}
  follow directly from  the following two known lemmas.

\begin{lem}\cite{So1}\label{lem-4-1}   Let $$f_n(x)=(x_1+ix_2)^n,\    \  x\in\sph.$$
Then   $f\in\mathcal{H}_n^d$ and
$$\|f_n\|_p\sim n^{-\ld /p},\    \   0< p<\infty.$$\end{lem}
\begin{lem}\cite[p.391]{Sz}\label{lem-4-2}
Let
$$g_n(x) =P_n^{(\f {d-3}2, \f {d-3}2)}(x\cdot e)$$
for a fixed point $e\in\sph$. Then
$g_n\in\mathcal{H}_n^d$, and
$$\|g_n\|_p \sim \begin{cases}
n^{\f {d-3}2} n^{-\f {d-1}p},\  & p> \f { 2(d-1)}{d-2},\\
n^{-\f12}(\log n)^{\f1p}, \   &p=\f { 2(d-1)}{d-2},\\
n^{-\f12},  & p<\f { 2(d-1)}{d-2}.\end{cases}$$
\end{lem}

For the proof of the upper estimates, we let
 $P\in\HH_n^d$.    The crucial tool in our proof is Theorem \ref{thm-2-3}, where we recall that $T_n P=P$ for all $P\in\HH_n^d$. We   consider the following cases:\\

{\it Case 1.}   \   \  $1\leq p\leq q\leq p'$.\\

In this case,  $1\leq p\leq 2\leq p'$, and the stated upper estimate for $q=p'$ follows directly from Theorem \ref{thm-2-3}.
In general, for $p\leq q \leq p'$, let $\t\in [0,1]$ be such that
$\f 1q=\f \t p+\f {1-\t}{p'}.$
Then
 by the log-convexity of the $L^p$-norm, we have
$$\|P\|_q\leq \|P\|_p^\t \|P\|_{p'}^{1-\t} \leq C n^{\l(\f 1p-\f1{p'})(1-\t)}\|P\|_p\leq C n^{\l( \f 1p-\f1q)}\|P\|_p,$$
which is as desired in this case.\\

{\it Case 2.}   \   \  $0<p\leq 1$ and $p<q$.\\

In this case, note that
$$\|P\|_1\leq \|P\|_p^p \|P\|_{\infty} ^{1-p}\leq
 C n^{\l(1-p)} \|P\|_p^p \|P\|_1^{1-p}.$$
It follows that
$$\|P\|_1\leq C n^{\l(\f 1p-1)}\|P\|_p,\    \  0<p\leq 1,$$
which, in turn,  implies that for  $p< q$ and  $\f 1q= \f {1-\t}p$,
    $$\|P\|_q \leq \|P\|_\infty^\t \|P\|_p^{1-\t}\leq
 C n^{\l \t } \|P\|_1^\t \|P\|_p^{1-\t}\leq C n^{\l (\f 1p-\f1q)}\|P\|_p.$$

 \vspace{5mm}

 {\it Case 3.}\  \  $1\leq  p\leq 2$ and  $q\ge p'$.\\

The desired estimate in this case follows directly
from  the first and the second parts of  Theorem \ref{thm-2-3} since $T_n P=P$ for all $P\in\HH_n^d$.\\

 {\it Case 4.}\  \ $2\leq p\leq  2+\f 1\l$ and  $q\ge 2+\f 2\l$.\\

 For $P\in\HH_n^d$, by the already proven cases it follows that
 \begin{align*}
    \|P\|_q&\leq C  n^{(d-2) (\f 12-\f1q)-\f1q}\|P\|_2\leq
    C  n^{(d-2) (\f 12-\f1q)-\f1q}\|P\|_p.
 \end{align*}

 {\it Case 5.}\  \  $2+\f 1\l<p<q\leq \infty$.\\

 The reverse H\"{o}lder inequality in this case follows directly from the corresponding Nikolskii inequality for spherical polynomials given by (\ref{Nikol}).\\

 {\it Case 6.}\  \
$d=3$  and  $2\leq p< 4=2+\f1\l$.\\

The  proof  in this case relies on the following result of Sogge \cite{So1}:

\begin{lem}\label{4-1-lem}If $d=3$, $\f 43<p<4$ and $q=3p'=p'(1+\f 1\l)$, then
$$\|\proj_n f\|_q \leq C n^{\f 12-\f 2q}\|f\|_p.$$
\end{lem}

Now we return to the proof in Case 6.
   Again, in view of Lemmas \ref{lem-4-1}  and \ref{lem-4-2},  it is enough to prove the upper estimates.
    Assume first that  $q\ge 3p'$.
Let $2\leq p< p_1 < 4$  and let $\t\in [0,1]$ be such that $$ \f 1p = \f {1-\t}{p_1}+\f \t 2.$$
Set $q_1=3p_1'$.
Then by Lemma \ref{4-1-lem},
\begin{equation}\label{4-1}
    \|Tf\|_{q_1} \leq C n ^{ \f 12 - \f 2 {q_1}}\|f\|_{p_1}.
\end{equation}
For $q\ge 3p'>3p_1'=q_1$, let $q_2\ge q$ be such that
$$\f 1q = \f {1-\t}{q_1} +\f \t{q_2}.$$
Then
\begin{align*}
\f13\ge \f 1{3p} +\f 1q =\t(\f 16 +\f 1{q_2}) + \f 13(1-\t)=\t (\f 1q_2-\f16)+\f13.
\end{align*}
This implies that $q_2\ge 6=2+\f 2\l$, hence by (ii) of Theorem \ref{thm-2-3},
\begin{align}\label{4-22}
\|T_n f\|_{q_2}\leq C n^{\f 12 -\f 2{q_2}} \|f\|_2.
\end{align}
Thus, using \eqref{4-1}, \eqref{4-22}, and the Riesz-Thorin theorem, we obtain
$$ \|T_n f\|_q \leq C n^{\f 12 -\f 2q}\|f\|_p,$$
which implies the desired estimate for the case of $q\ge 3p'$.

The case of $p<q< 3p'$ can be treated similarly. In fact, let $p_1, q_1$ and $\t$ be as above. Observing  that $\f12 -\f 2 {q_1} =\f 12 (\f 1{p_1}-\f 1{q_1})$, we may rewrite \eqref{4-1} as
\begin{equation*}
    \|Tf\|_{q_1} \leq C n ^{ \f 12 (\f 1{p_1}-\f 1{q_1})}\|f\|_{p_1}.
\end{equation*}
Furthermore, we may choose $p_1>p$ to be very close to $p$ so that $q < q_1=3p_1'<3p'$. Let $q_3\leq q $ be such that
$$\f 1q = \f {1-\t}{q_1} +\f \t{q_3}.$$
Then
\begin{align*}
\f13< \f 1{3p} +\f 1q =\t(\f 16 +\f 1{q_3}) + \f 13(1-\t)=\t (\f 1q_3-\f16)+\f13.
\end{align*}
Hence $2<q_3 <6$, and using
(i) of Theorem \ref{thm-2-3},  we deduce
\begin{align*}
\|T_n f\|_{q_3}\leq C n^{\f 12 ( \f 12 -\f 1{q_3})} \|f\|_2.
\end{align*}
The stated estimate for $p<q<3p'$ then follows by the Riesz-Thorin interpolation theorem.
\hb

\section{Applications: Fourier inequalities}

\subsection{ The restriction conjecture. } One of the most challenging  problems  in  classical Fourier analysis is
the restriction conjecture, which states that if $1 \le p < \frac{2d}{d+1}$ and $q \leq  \frac{ d-1}{
d+1} p'$,  then  there exists a constant $C$ depending only on $p, q, d$ such that \begin{equation}\label{rc0}
\frac{\| \widehat{F}	\|_{L^q(\mathbb{S}^{d-1})}}
{\| {F}	\|_{L^p(\mathbb{R}^{d})}}
\le C,\      \    \  \forall F \in C^\infty_0(\mathbb{R}^d),
 \end{equation}
 where  $\hat{F}(\xi): =\int_{\RR^d} F(x) e^{-2\pi i x\cdot \xi}\, dx,$  $  \xi\in\RR^d.$
    This conjecture has been  completely  proved  only  in the case of $d=2$.  We refer to  the book \cite[Chapter IX]{stein} for more background information of this problem.

 De Carli and Grafakos \cite{grafakos}  recently proved that  the restriction conjecture is valid for all  functions  $F$ that can be expressed in the form
  $$F(x)=f(\|x\|) \|x\|^n  g_n\bl(\f{x}{\|x\|}\br),\   \   \  n=0,1,\cdots$$ with  $f(\|\cdot\|)\in C_0^\infty(\RR^d) $ and $g_n\in\HH_n^d$ being  given in   \eqref{1-3} . Using \thmref{cor2} (i), and following  the argument of \cite{grafakos}, we   may   conclude here  that
   the restriction conjecture  holds for a wider class of functions
   $$F\in \bigcup_{n=0}^\infty \Bl\{ f(\|x\|) \|x\|^n  Y_n\bl(\f{x}{\|x\|}\br):\   \
   f(\|\cdot\|)\in C_0^\infty(\RR^d),\   \   Y_n\in\HH_n^d\Br\}.$$
  Indeed, it was shown in \cite{grafakos}  that for $F(x)=f(\|x\|)\|x\|^n
  Y_n(x/\|x\|)$ with $f\in C_0^\infty (\RR^d)$ and $Y_n \in\HH_n^d$,
\begin{align}
\frac{\| \widehat{F}	\|_{L^q(\mathbb{S}^{d-1})}}
{\| {F}	\|_{L^p(\mathbb{R}^{d})}}
&=
\frac
{
\Big|
\int_0^{\infty} f (r) J_{\frac d2-1+n}(r)r^{\frac d2+n} dr\Big|
}
{
\Big(
\int_0^{\infty} |f(r)|^p r^{d-1+np} dr\Big)^{1/p}
}
\frac{ \|Y_n\|_{L^q(\mathbb{S}^{d-1})}  }{ \|Y_n\|_{L^p(\mathbb{S}^{d-1})}  }
\nonumber
\\
&\leq
C n^{(d-1)(\f 1 2-\f 1 p)+\f 1 {p'}}\frac{ \|Y_n\|_{L^q(\mathbb{S}^{d-1})}  }{ \|Y_n\|_{L^p(\mathbb{S}^{d-1})}  }, \label{5-2}\end{align}
  where $J_n(r)$ is the Bessel function of the first kind.
However,   according to (i) of \thmref{cor2} ,  we obtain that for $1 \le p < \frac{2d}{d+1}$ and $q \leq \frac{ d-1}{
d+1} p'$,
\begin{align*}
\text{RHS of \eqref{5-2}} \leq C
\sup\limits_{m\ge 1}
m^{(d-1)(\frac12-\frac1p)+\frac{1}{p'}+\frac{d-2}{2} (\f {1}{p}-\f 1q)} \leq C.
\end{align*}

\subsection{The sharp Pitt inequality}
The following  sharp  Pitt inequality has been recently proved in \cite{DeC}:

\begin{thm}\label{T-Pitt2} If $1\leq p\leq 2$ and  $s=
(d-1)\left(\frac 12-\frac 1p\right)$,  then for every $Y_k\in \CH_k^d$ and  every radial $f\in {\mathcal S}(\R^d)$,
 the Pitt  inequality
\begin{equation}\label{e-newpitt}
\|\,|y|^{-s} \wh{f Y_k} \|_{L^{p'} (\R^d)} \leq C \|\,|x|^s fY_k \|_{L^p(\R^d)}
\end{equation}
holds with  the   best  constant
\begin{equation}\label{e-B}
C=(2\pi)^{\frac d2} 2^{\frac{1}{2}-\frac{1}{p'}} \frac{p^{\frac{(2 k+d-1)p+2}{4p}} \Gamma \left(\frac{(2 k+d-1)
p'+2}{4}\right)^{\frac{1}{p'}}}{(p')^{\frac{(2 k+d-1)p'+2}{4p'}} \Gamma
\left(\frac{(2 k+d-1) p+2}{4}\right)^{\frac 1p}}\sup_{Y_k\in \CH_k^d} \frac{\|Y_k\|_{L^{p'}(\SS^{n-1})}}{\|Y_k\|_{L^ p(\SS^{n-1})}}.
\end{equation}
\end{thm}

According to  \thmref{cor2}, we have
$$\sup_{Y_k\in \CH_k^d} \frac{\|Y_k\|_{L^{p'}(\SS^{n-1})}}{\|Y_k\|_{L^ p(\SS^{n-1})}}\sim
k^{(d-2)(\frac{1}{p}-\frac{1}{2})},$$
whereas  only the weaker estimate \eqref{1-2} was obtained in  \cite{DeC}.

\end{document}